\input amstex
\documentstyle{amsppt}
\input bull-ppt
\keyedby{harvey/amh}   
\define\a{\alpha}
\predefine\barunder{\b}
\redefine\b{\beta}

\predefine\dotunder{\d}
\redefine\d{\delta}
\predefine\dotaccent{\D}
\redefine\D{\Delta}

\define\equdef{\overset\text{def}\to=}
\define\f{\frac}

\define\G{\Gamma}

\define\lb\{{\left\{}
\define\la{\lambda}
\define\La{\Lambda}

\define\lra{\longrightarrow}
\define\notS{\hbox{$S\negthickspace\negthickspace /\;$}\!\!}
\define\nots{\hbox{$\!{\ssize 
S\negthickspace\negthickspace /\;}$}\!\!}

\define\ola{\overleftarrow}
\define\Om{\Omega}
\define\om{\omega}
\define\op{\oplus}
\define\oper{\operatorname}

\define\ora{\overrightarrow}
\define\ov{\overline}

\define\ox{\otimes}

\define\s{\sigma}

\define\th{\theta}

\define\wh{\widehat}
\define\wt{\widetilde}

\define\lp({\left(}
\define\rtp){\right)}
\define\slantline#1#2#3#4#5{\hbox to 0pt{
}}


\predefine\preSL{\SL}
\redefine\SL{{\Cal L}}

\define\Todd{\text{\bf Todd}}
\define\BA{{\bold A}}

\define\BC{{\bold C}}

\define\BH{{\bold H}}

\define\BP{{\bold P}}

\define\BR{{\bold R}}

\define\BZ{{\bold Z}}


\define\bbc{{\Bbb C}}

\define\rarr{\rightarrow}
\define\odd{\text{odd}}
\define\even{\text{even}}
\define\loc{L^1_{\text{loc}}}

\define\rank{\oper{rank}}
\define\Res{\oper{Res}}
\define\Hom{\oper{Hom}}
\define\Div{\oper{Div}}


\predefine\prech{\ch}
\redefine\ch{\text{ch}}

\define\tr{\text{tr}}

\predefine\prebh{\bh}
\redefine\bh{\bold h}

\define\bF{\bold F}

\define\btau{\boldtau}


\define\bvf{\bold\phi}


\define\arr{\longrightarrow}
\define\vf{\phi}
\define\9{\roman{TP}_s}
\define\TP{\roman{TP}_\infty}

\ratitle
\topmatter 
\cvol{31}
\cvolyear{1994}
\cmonth{July}
\cyear{1994}
\cvolno{1}
\cpgs{54-63}

\title A theory of characteristic currents\\ 
associated with a singular connection\endtitle
\shorttitle{Theory of characteristic currents}
\author Reese Harvey and H. Blaine Lawson, Jr.\endauthor 
\shortauthor{Reese Harvey and H. B. Lawson, Jr.}
\address Department of Mathematics, Rice University,
Houston, Texas 77251-1892\endaddress
\address Department of Mathematics, SUNY at Stony Brook,
Stony Brook, New York 11794-3651\endaddress
\thanks Research of both authors was supported by the
National Science Foundation\endthanks
\date October 6, 1992 and, in revised form, August 20, 
1993\enddate
\abstract This note announces a general construction of
characteristic currents for singular connections on a 
vector bundle.
It develops, in particular, a Chern-Weil-Simons  theory 
for smooth bundle
maps $\alpha : E \rightarrow F$ which, for smooth 
connections on $E$ and
$F$, establishes formulas of the type
$$
\phi \ = \ \text{\rm Res}_{\phi}\Sigma_{\alpha} + dT.
$$
Here $\phi$ is a standard charactersitic form, 
$\text{Res}_{\phi}$
is an associated smooth ``residue'' form computed 
canonically in
terms of curvature, $\Sigma_{\alpha}$ is a  rectifiable 
current
depending only on the singular structure of $\alpha$, and 
$T$ is a
canonical, functorial transgression form with coefficients 
in
$\loc$. The theory encompasses such classical topics as:
Poincar\'e-Lelong Theory, Bott-Chern Theory, Chern-Weil 
Theory, and
formulas of Hopf.  Applications  include:\ \ 
a new proof of the Riemann-Roch Theorem for vector bundles 
over algebraic
curves,  a $C^{\infty}$-generalization of the
Poincar\'e-Lelong Formula,  universal formulas for the 
Thom class
as an equivariant characteristic form (i.e., canonical 
formulas for
a de Rham representative of the Thom class of a bundle with
connection), and  a Differentiable 
Riemann-Roch-Grothendieck Theorem at
the level of forms and currents. A variety of formulas 
relating
geometry and characteristic classes are deduced  as direct
consequences of the theory.\endabstract

\endtopmatter

\document

\heading{1. Introduction}\endheading 
Our purpose is to announce a
generalized Chern-Weil theory for singular connections on 
a vector
bundle, particularly those which arise from bundle maps.
Specifically, if $E$ and $F$ are smooth vector bundles with
connection over a manifold $X$, the theory associates to 
each
homomorphism $\a : E\arr F$ a family of $d$-closed 
characteristic
currents on $X$ defined canonically in terms of the 
curvature of the
bundles and the singularities of the map $\a$. 

One can focus attention either on $E$ or on $F$ (and 
retrieve, when
$\a\equiv 0$, the standard theory for $E$ or $F$). Suppose 
the focus
is on $F$, and fix a {\it characteristic polynomial} 
$\vf$, i.e., an
Ad-invariant polynomial on the Lie algebra of the 
structure group
of $F$. The new  theory associates to $\vf$ a $d$-closed 
current
$\vf(\Omega_{F,\a})$ which is defined in terms of 
curvature and the
singular structure of $\a$ and which represents the
characteristic class $\bvf(F)$. Moreover, the theory 
constructs a
canonical and functorial transgression current $T = 
T(\vf,\a)$ with
the property that  
$$
\vf(\Omega_{F}) - \vf(\Omega_{F,\a}) \, = \, dT.
$$
 
When $\rank(E) = \rank(F)$, this often yields an
equation of the form 
$$
\vf(\Omega_{F}) - \vf(\Omega_{E}) \ =\ \Res_{\vf}[\Sigma] +
 dT
$$
where $[\Sigma]$ is the current given by integration over 
the singular
set $\Sigma\equdef \{x\in X\, : \ \a_x \text{\ is not an
isomorphism}\}$ and where $\Res_{\vf}$ is a smooth 
``residue'' form
which is expressed canonically in terms of the curvatures
$\Omega_{E}$ and $\Omega_{F}$ of $E$ and $F$.
One can think of $S_{\vf} =\Res_{\vf}[\Sigma]$ as the 
characteristic 
current associated to $\a : E\arr F$ which represents the 
class
$\bvf({F}) - \bvf({E})\in H^*_{\text{de Rham}}(X)$.

When rank($E$) $>$ rank($F$), there are analogous formulas:
$$
\vf(\Omega_E) - \vf(\Omega_0) \ = \   \Res_{\vf}[\Sigma] + 
dT
$$
where $ \vf(\Omega_0)$ is a closed differential form with
$\loc(X)$-coefficients whose cohomology class is more subtly
determined  by the global geometry.

\heading{2.  An example: Complex line bundles}\endheading
Suppose that $E$ and $F$ are complex
line bundles over an oriented manifold $X$ and that $\a : 
E\arr F$ is a
smooth bundle map. If $\a$ vanishes nondegenerately, we 
define its 
divisor to be the current $\Div(\a) = [\Sigma]$  
associated to the
oriented codimension-2 submanifold $\Sigma=\{x\in X\, : \, 
\a_x=0\}$. 
Set  
$$
f = \tfrac{i}{2\pi}\Omega_F \quad\text{and}\quad
e = \tfrac{i}{2\pi}\Omega_E
$$
where $\Omega_E$ and $\Omega_F$ are the curvature 2-forms 
of the
connections. Then it is shown that there exists a 1-form 
$\s$ on $X$
with $\loc(X)$-coefficients such that
$$
f - e  \ =\ \Div(\a) + {\tsize\frac{i}{2\pi}} d\s,\tag1
$$
where $\s$ is defined by $D\a=\s\a$. Passing to 
cohomology, we capture
the first Chern class  $c_1(F) - c_1(E)  = [f]-[e] = 
[\Div(\a)]\in
H^2_{\text{de Rham}}(X)$. Note that when $E$ is trivial, 
equation
(1) represents a $C^{\infty}$ enhancement of the classical
Poincar\'e-Lelong Formula for sections of the bundle $F$.

More precisely, formula (1) is shown to hold for any 
bundle map $\a$
which is {\it atomic}, that is, for which
$$
\f{da}{a}\in\loc(X),
$$
where $a$ is the  complex-valued function which represents 
$\a$ in
the given local frames. This condition is independent of 
the choice
of frames, as is the closed current $\Div(\a) = 
d\lp(\f1{2\pi i}\f{da}a\rtp)$ (cf\. [HS]).  Note that 
locally 
$\s = \f{da}a+\om_F-\om_E$.
 The following general result
is formula (1) when $\vf(u) = u$.

\proclaim{Theorem A} Fix a  polynomial $\vf(u)\in\BC[u]$ 
in one
indeterminate. Then for any atomic bundle map $\a : E \arr 
F$ as above,
there exists a differential form $T$ with 
$\loc$-coefficients on $X$
such that 
$$
\vf(f) - \vf(e) \ = \ \left\{\frac{\vf(f) -
\vf(e)}{f-e}\right\}\Div(\a) + dT\tag2
$$
where $\frac{\vf(f) - \vf(e)}{f-e}$ is the obvious 
polynomial in $e$ and
$f$. In fact, the $\loc(X)$-form $T$ is given explicitly by
$$ 
T = {\tsize\frac{i}{2\pi}}
\left\{\frac{\vf(f) - \vf(e)}{f-e}\right\}\s.
$$
\endproclaim

Combining this result with the kernel-calculus of Harvey 
and Polking
\cite{HP} yields a new proof of the Riemann-Roch Theorem 
for vector
bundles over algebraic curves.

\heading{3. The general procedure; approximation  
modes}\endheading
Let $E$ and $F$ be smooth vector bundles with connections 
$D_E$ and
$D_F$ over a manifold $X$ where, for simplicity in the 
following
discussion, we assume $\rank(E)\le\rank(F)$. Now for 
bundle maps
$$
\a : E\arr F\quad\text{and}\quad\b : F\arr E
$$
these connections can be transplanted via the formulas
$$
\ola D^{\a,\b} = \b{\circ}D_F{\circ} \a \ +\ (1-\b\a)
{\circ} D_E \quad \text{and} \quad 
\ora D^{\a,\b} =  \a{\circ} D_E{\circ} \b \ +\ 
D_F{\circ}(1-\a\b)
$$
to define induced connections on $E$ and $F$ respectively. 
When $E$ and
$F$ are isomorphic and $\b = \a^{-1}$, one recovers the 
standard
gauge transformations. If $\a$ is injective, one can 
introduce
metrics on $E$ and $F$ and define 
$$
\b = (\a^*\a)^{-1} \a^*,
$$
yielding a {\it pull back connection} $\ola D$ and a
{\it push forward connection} $\ora D$. In general this 
procedure breaks down on
the {\it singular set} $\Sigma\equiv\{x\ :\ \a_x\ \text{is 
not
injective}\}$, since $\b$ becomes singular on $\Sigma$. To 
remedy this, we
choose an {\it approximation mode} by fixing a 
$C^{\infty}$-function
$\chi: [0,\infty]\arr [0,1]$ with $\chi'\geq 0$, $\chi(0) 
= 0$, and
$\chi(\infty) = 1$. We then define a smooth approximation 
$\b_s$ to $\b$
by 
$$
\b_s = \chi\left(\frac{\a^*\a}{s}\right)\b
\quad \text{for}\ s>0.
$$
Plugging $\b_s$  into the formulas above gives  a family 
of smooth
connections $\ora D_s$ on $F$ (and $\ola D_s$ on $E$). As
$s\rightarrow 0$, $\b_s\rightarrow\b$ uniformly on compacta 
in $X-\Sigma$. If $\chi(t)\equiv 1$ for $t\geq 1$, then 
$\b_s = \b$
outside the neighborhood $\{x\in X : \|\a_x\|^2 <s \}$ of 
$\Sigma$.
Such a choice of $\chi$ is called a
{\it compactly supported approximation mode}. Another 
important case
is the {\it algebraic approximation mode} where $\chi(t) =
t/(1+t)$.  Here the family of connections $\ora D_s$ is 
directly
related to the {\it Grassmann graph construction} of
R\.~MacPherson. When working in this mode with the 
tautological
bundle map $\a:\pi^*E\arr\pi^*F$ over the total space of
$\pi:\Hom(E,F)\rightarrow X$, the family $\ora D_s$ 
extends smoothly
to a  fibrewise compactifiation of $\Hom(E,F)$. This has 
strong
consequences for characteristic forms in the curvature of 
$\ora D_s$.

Let $\ora\Omega_s$ denote the curvature of $\ora D_s$, and 
fix an
Ad-invariant polynomial $\vf$ on the  Lie algebra of the 
structure
group of $F$.

\subheading{Definition A} Suppose that the limit 
$$
\vf(\!(\ora D)\!)\equdef\underset{s\rarr
0}\to{\lim}\vf(\ora\Omega_s)
$$
exists as a current on $X$. Then $\vf(\!(\ora D)\!)$ is 
called the 
{\it $\bvf$-characteristic current associated to the 
singular
push forward connection $\ora D$ on $F$}. An analogous 
definition
holds for the singular pull back connection $\ola D$ on $E$.

Note that $\vf(\!(\ora D)\!)$ is automatically $d$-closed 
and
represents the $\vf$-characteristic class of 
$F$. Over the subset
$X-\Sigma$, $\ora D$ is a smooth connection, and 
$\vf(\!(\ora D)\!)$
is a smooth differential form which we denote by
$\vf(\ora\Omega_0)$. Note that when $\rank(E)\ = 
\rank(F)$, \ $\ora
D$ is gauge equivalent via $\a$ to $D_E$ over $X-\Sigma$, 
and
therefore $\vf(\ora\Omega_0) = \vf( \Omega_E)$  extends 
smoothly
across the singular set $\Sigma$. When $\rank(E)\ < 
\rank(F)$, 
[HL1] establishes conditions on $\a$ which guarantee that
$\vf(\ora\Omega_0)$ extends across  $\Sigma$ as a $d$-closed
$\loc$-form. This gives a decomposition 
$$ 
\vf(\!(\ora D)\!) = \vf(\ora\Omega_0) + S_{\vf}\tag3
$$
where $\vf(\ora\Omega_0) $ denotes the $\loc(X)$-extension 
and
where $S_{\vf} $ is a current on $X$ with the property that 
$$
d S_{\vf}  = 0 \quad\text{and}\quad
\text{supp}(S_{\vf} ) \subseteq \Sigma.
$$
Each term in (3) represents a de Rham cohomology class on 
$X$ which
can be nonzero even when $E$ and $F$ are trivial bundles. 

The detailed structure of $S_{\vf} $ and its independence of
approximation mode is established by considering a family
of transgression forms $T_s$ with 
$$
dT_s = \vf(\Omega_F) - \vf(\ora\Omega_s)
$$
and finding conditions on $\a$ so that
$T\equdef\underset{s\rarr 0}\to{\lim} T_s$ exists in
$\loc(X)$. The {\it transgression} $T$ satisfies
$$
dT\ =\ \vf(\Omega_F) - \vf(\ora\Omega_0) - S_{\vf}
$$
and is   functorial under appropriately transverse maps 
between
manifolds. In many cases $S_{\vf} $ can be written as
$$
S_{\vf} = \Res_{\phi}[\Sigma]
$$
where $\Res_{\phi}$ is a smooth form on $X$, expressed as 
a universal
Ad-invariant polynomial in $\Omega_E$ and $\Omega_F$. In 
particular 
$\Res_{\phi}$ is completely determined by computing its 
associated 
cohomology class in the universal setting.

Combining this gives a {\it Chern-Weil Theorem for bundle 
maps}:
$$
\phi(\Omega_F) - \phi(\ora\Omega_0) \ =\  
\Res_{\phi}[\Sigma] + dT 
\tag4
$$
where $[\Sigma]$ is a current canonically determined by 
the singular
structure of $\a$. Special cases are  discussed in the 
following
sections.

\heading{4. Universal Thom classes}\endheading 
When $E$ is the trivial line
bundle, $\a$ becomes a cross-section of $F$.  

\subheading{Definition B} The cross-section $\a$ is said 
to be 
{\it atomic} if whenever   $\a$ is written locally as
$\a = \sum a_je_j$ in terms of a local frame field $e$, the
$\BR^m$-valued function $a = (a_1,\dots,a_m)$ satisfies 
$$
\frac{da^I}{|a|^{|I|}}\in\loc(X) \quad \text{for all \ 
}|I| < m.
$$ 

Elementary criteria for atomicity are given in \cite{HS}. 
Roughly
speaking, any smooth  section which vanishes 
``algebraically'' on a
set of the appropriate ``Minkowski codimension'' is 
atomic. For
example, any real analytic section with zeros of 
codimension $\geq 
m$ is atomic. 
 
It is proved in \cite{HS} that the vanishing of an atomic 
section
$\a$ determines a unique, $d$-closed current $\Div(\a)$ of 
real
codimension-$m$, called the {\it divisor} of $\a$. 
Locally, the
divisor is defined by 
$$ 
\Div(\a)\ =\ d(a^*\th)
$$
where $\th$ is the normalized solid angle kernel on 
$\BR^m$. This
current $\Div(\a)$ is integrally flat, and in particular, 
when its mass
is finite, it is a rectifiable cycle in the sense of 
Federer \cite{F}.
 
Suppose that $F$ is either complex, or real and oriented, 
and that
$\vf$ is the top Chern polynomial $\wt{\det}(A)=\det(\f 
i{2\pi}A)$
or the Euler polynomial
$\wt{\text{Pf}}(A)=\text{Pfaff}\lp(-\f1{2\pi}A\rtp)$ 
respectively. (If
$\dim_{\BR}F$ is odd, then $\vf \equiv 0$.)  When $\a$ is 
atomic, 
the theory produces a canonical $\loc(X)$-form $\s$ on $X$ 
called
the {\it spherical potential} such that 
$$ 
\vf(\Omega_F)\,-\,\Div(\a)\ =\ d\s. 
$$
Furthermore, for each approximation mode there is a smooth
family of connections $\ola D_s$, $0<s\leq\infty$ on $F$, 
so that  
$\tau_s\equiv\vf(\ola \Omega_{s})$ satisfies 
$$
\tau_{\infty} = \vf(\Omega_F)
\quad\text{and}\quad
\lim_{s\rarr 0}\tau_s\ =\ \Div(\a).
$$
There is also a family of $\loc$-forms $\s_s$ with 
$\underset{s\rarr
0}\to{\lim}\s_s = \s$ in $\loc(X)$ and
$$
\tau_s\,-\,\Div(\a) \ =\ d(\s-\s_s).
$$

It is useful to view this construction on the
total space of the bundle $\pi: F\arr X$.  Consider the 
pull back
$\bF = \pi^*F$ with the pull back connection. Over $F$ 
there is a
{\it tautological cross section} $\boldalpha$ of $\bF$ given
by $\boldalpha(v) = v$. This section is atomic,  and the 
theory
applies. Each approximation mode $\chi$ gives a smooth 
family of
closed differential forms $\btau_{\bold s}$, $0 < 
s\leq\infty$ on $F$,
which are expressed canonically in terms of the 
connection, and
which represent the {\it Thom class} of $F$. For example, 
when $F$
is real of dimension $2n$ and $\chi(t) = 1-1/\sqrt{1+t}$, 
$\btau_{\bold s}$ is given by the formula
$$
\btau_{\bold s}\ =\frac{s}{\sqrt{|u|^2 +s^2}}
\wt{\text{Pf}}\!\left({\frac{Du^t Du}{\ssize |u|^2+s^2}}
- \Omega_F \right)
$$
where $u$ represents  the tautological section and $Du$ is 
its
covariant derivative.

\proclaim{Theorem B} For any approximation mode the form 
$\btau_{\bold s}$
is a closed $2n$-form on $F$ which dies at
infinity, is integrable on the fibres with integral one, 
restricts
to be $\wt{{\roman{Pf}}}(\Om_F)$ on $X$, converges as 
$s\rarr 0$ to
the current $[X]$ represented by the zero-section, and 
converges
as $s\to\infty$ to $\wt{{\roman{Pf}}}(\Om_F)$.
\endproclaim

The {\it Thom form} $\boldtau_{\bold s}$ can be written as 
the Chern-Weil image
of a universal form  in the equivariant de Rham complex of 
$\BR^m$.
In each   approximation mode there is  an explicit universal
formula for the Thom forms $\boldtau_{\bold s}$. If 
$\chi(t)\equiv 1
\ \text{for\ } t\geq 1$, we have  supp$(\btau_{\bold s}) 
\subset\{v\in
F\,:\,|v|\leq s\}$. In particular, $\btau_{\bold s}$ has 
{\it compact support in each fibre.} 

There are analogous Thom forms associated to the top Chern 
form when 
$F$ is complex.

\heading{5. Rectifiable 
Grothendieck-Riemann-Roch}\endheading 
Another  
basic formula arises from the theory when considering an
atomic section $\a\in\G(V)$ of an even-dimensional, real 
vector bundle
$V\arr X$ with spin structure. Clifford multiplication by 
$\a$
determines a bundle map $\alpha:
{\notS}^+\arr{\notS}^-$ between the positive
and negative complex spinor bundles of $V$.
Consider the function on matrices $\text{\bf ch}(A)\equdef
\text{trace}\left\{\exp(\frac{1}{2\pi i} A)\right\}$ which 
gives the
Chern character. Suppose ${\notS}^+$ and ${\notS}^-$ carry 
connections induced
from a metric connection on $V$, and let 
$\Omega_{{\nots}^\pm}$,
$\Omega_{V}$ denote the curvature matrices of these 
connections.

\proclaim{Theorem C} The following equation of forms and 
currents 
holds on $X$\RM:
$$
\text{\bf  ch}\!\left(\Omega_{{\nots}^+}\right)
- \text{\bf ch}\!\left(\Omega_{{\nots}^-}\right) \ =\
\wh\BA(\Omega_{V})^{-1}\Div(\a) + dT\tag5
$$
where
$$
\wh\BA(\Omega_{V})^{-1}\ =\ {\det}^{\f 1 
2}\!\left\{\frac{\sinh(
{\ssize\frac{1}{4\pi}} \Omega_{V})}
{{\ssize\frac{1}{4\pi}}\Omega_{V}}\right\}
$$
is the series of differential forms on $X$ which canonically
represent, via Chern-Weil Theory, the inverse 
$\wh\BA$-class of $V$,
and $T=\wh \BA^{-1}(\Om_V)\s$ where $\s$ is the spherical 
potential.
More generally, for any complex bundle $E$ with connection 
over $X$
and $T=\roman{ch}\,(\Om_E)\wh \BA^{-1}(\Om_V)\s$ one 
has\,\RM:
$$
\text{\bf ch}\!\left(\Omega_{{\nots}^+\otimes E}\right)
- \text{\bf ch}\!\left(\Omega_{{\nots}^-\otimes E}\right) 
\ =\
\text{\bf ch}(\Omega_E)\wh\BA(\Omega_{V})^{-1}\Div(\a) + 
dT.\tag6
$$
\endproclaim

Equation (6) raises the Differentiable Riemann-Roch Theorem
for embeddings \cite{AH} to {\it  the level of 
differential forms} 
and extends it to oriented subcomplexes which arise as 
divisors of
some cross section of a bundle. One also obtains 
approximating families
as with the Thom forms above.

There are formulas analogous to (6) when $V$ is complex or
$\text{Spin}^c$. The complex case  gives a canonical 
version of
a classical Grothendieck Theorem at the level of forms
and currents. To be specific, let $j:Y \hookrightarrow X$ 
be a
proper almost complex embedding of almost complex 
manifolds with
normal bundle $N$. If the tangent bundles $TY$ and $TX$ 
are given
connections compatible with the complex structures,
then one has the following equation of forms and currents 
on $X$:
$$\align
& \left\{\text{\bf ch}
\!\left(\Omega_{\left(\La^{\even}N^*\right)\otimes 
E}\right)   
- \text{\bf 
ch}\!\left(\Omega_{\left(\La^{\odd}N^*\right)\otimes
E}\right)\right\}  \wedge\Todd \left(\Omega_{TX}\right)\\
&\qquad= \text{\bf  ch}(\Omega_E)\wedge\Todd(\Omega_{TY}) 
[Y]  + dT\endalign
$$
for any vector bundle with connection $E$ over $Y$.  By 
passing to
cohomology, this formula yields the commutativity of the 
diagram
$$
\CD
K(Y) @>j_{\,!}>> K(X)  \\
@V{\text{\bf  ch}(\,\cdot\,)\wedge \Todd(Y)}VV     
@VV{\text{\bf  ch}(\,\cdot\,)\wedge
\Todd(X)}V    
 \\ H^*(Y) @>>j_{\,!}> H^*(X)
\endCD
$$
where the ${j_!}$ represent the Gysin ``wrong way'' maps 
in K-theory and
cohomology.

In  these special Clifford multiplication cases this 
formalism has
some similarities with Quillen's calculus of 
superconnections
\cite{Q} as developed in \cite{MQ, BV, BGS*} and
elsewhere. However, even in these cases there are 
substantial
differences. We are concerned with transgressions and  
convergence 
questions under the
weak atomic hypothesis and with the explicit structure of 
the
limiting currents. Our theory also allows for a quite 
general choice
of approximation mode.  Choosing $\chi(t) = 1-e^t$ puts us 
closest to
Quillen's theory. However, choosing $\chi$ with 
$\chi(t)\equiv 1$
for $t\geq 1$ gives approximations supported in small 
neighborhoods
of the singular set.  Approximations where $\chi(t) = t/(1+
t)$ 
admit nice compactifications (see \cite{Z1}) and are  
related to
MacPherson's  Grassmann graph construction.

\heading{6. Quaternionic line bundles}\endheading Let $\a: 
E\arr F$  be an
endomorphism of smooth quaternionic line bundles. The
classifying space for quaternionic line bundles is the 
infinite
quaternionic projective space $\BP^{\infty}(\BH)$ whose 
cohomology is a
polynomial ring on one generator 
$u\in H^2(\BP^{\infty}(\BH);\, \BZ)$ called the
{\it instanton class}.

\proclaim{Theorem D} Suppose $E$ and $F$ are provided with
connections which are compatible with the quaternionic 
structure, and
assume that $\a$ is atomic. Then for each $\vf\in\BR[u]$ 
there exists a
canonical $\loc(X)$-form $T$ with the property that 
$$
\vf(f) - \vf(e) \ = \ \left\{\frac{\vf(f) -
\vf(e)}{f-e}\right\} \Div(\a) + dT
$$
where
$$
f = \tfrac{1}{16\pi^2 }\tr\left\{\Omega_F^2\right\}
\quad\text{and} \quad  e =\tfrac{1}{16\pi^2}\tr
\left\{\Omega_E^2\right\}
$$
are the canonical representatives of the instanton class  
of $E$ and
$F$. 
\endproclaim

\heading{7. Degeneracy currents}\endheading  For each 
bundle map $\a:E \arr F$
and integer $k$ with $k<m = \rank(E) \leq n = \rank(F)$, we
define a current which measures the locus where 
rank$(\a)\leq k$. 
Specifically, let $\xi:G_r(E)\arr X$ be the Grassmann
bundle of $r$-planes in $E$  with $r=m-k$, and let $U\arr
G_r(E)$ be the tautological $r$-plane bundle. Note that 
$U$ is a
subbundle of $\xi^*E$ and that $\a$ pulls back to a bundle 
map
$\xi^*\a: \xi^*E\to\xi^*F$. We say that the bundle map 
$\a$  is
$k$-{\it atomic} if the restriction $\hat{\a} = 
\xi^*\a\big|_U$ is an
atomic section of $\roman{Hom}(U,\xi^*F)$ over $G_r(E)$.  
Under this
hypothesis we define the $k^{\text{th}}$ {\it degeneracy
current} of $\a$ to be the current push-forward
$$
{\Bbb D}_k(\a) \,=\, \xi_* \Div(\hat{\a}).
$$
When $\a$ is algebraic, $ {\Bbb D}_k(\a)$ is a cycle whose 
class in the
Chow ring $A^*(X)$ is the one defined by Fulton in 
\cite{Fu}. For
given connections on $E$ and $F$ there are canonical 
families of
smooth forms $\9$ and $\loc$-forms $S_s$ on $X$, for
$0<s\leq \infty$, such that
$$
{\Bbb D}_k(\a) \,=\, \9 + dS_s\tag7
$$
where $\underset{s\rarr 0}\to{\lim} S_s = 0$ in $\loc$ 
and where in the complex case
$$
\TP\,=\, \text{det}_{m\times m}
(( c(E-F)_{n-i+j})).
$$
Here $c(E-F)\equdef 1 +c_1 +c_2+\dots = 
c(\Omega_F)/c(\Omega_E)$ in
the ring of even forms on $X$.  There is an analogous 
formula in the
real case. These {\it Thom-Porteous families} $\9$ are 
established
in  \cite{HL2}.

\heading{8. Geometric formulas}\endheading The results 
above give a wide
variety of formulas relating characteristic forms to 
singularities
of maps.  For example, consider  $\a : \underline{\BC}^{k+
1}\rightarrow F$
corresponding to $(k+1)$-sections $s_0,\dots ,s_k$ of a
complex bundle with connection $F$.  If $\a$ is 
$k$-atomic, the top
degeneracy current $ {\Bbb D}_k(\a) 
\equdef {\Bbb L}{\Bbb D}(s_0,\dots ,s_k)$ is defined and 
measures the
locus of linear dependence of these sections.  From (7) we 
have
$$
c_{n-k}(\Omega_F) - {\Bbb L}{\Bbb D}(s_0,\dots ,s_k) \ =\ dT
$$
where $T$ is a canonical $\loc$-form on $X$. There are 
corresponding
formulas for higher-order degeneracies. 

A fundamental case occurs when $\a=df: TX\lra f^*TY$ is 
the differential
of a smooth map $f: X\lra Y$ between manifolds. This 
yields, for
example, a $C^{\infty}$-version of the formula for the
global Milnor current (cf. [Fu, 14.1.5]). When $X$ and $Y$ 
are
oriented 4-manifolds, it  yields a formula relating 
$f^*p_1(Y)-p_1(X)$
to a weighted sum of isolated ``singular points'' where 
rank($df$) =
2  (cf\. \cite{M2} and \cite{R}). There are also formulas 
of type
(7) explicitly relating Pontrjagin forms to the  
singularities of
projections $X^n \subset {\bold R}^{n+n'} \arr {\bold 
R}^{k}$.

Applications to CR-geometry include the following.  
Consider a
generic 
smooth map $f:M\rightarrow\BC^{k+1}$ of an oriented 
riemannian
$n$-manifold $M$ where $n-k = 2\ell >0$. Then there are  
formulas
$$
p_{\ell}(\Omega_M) = (-1)^{\ell}{\Bbb C}r (f) + dT
$$
where $p_{\ell}(\Omega_M) $ is the $\ell^{\text{th}}$ 
Pontrjagin form of  $M$
and where ${\Bbb C}r (f)$ is a current associated to the 
{\it CR
critical set}:\ \ $\{x\in M \,:\, df : 
TM\otimes_{\BR}\BC\arr\BC^{k+1}$
 is not surjective\}. For example, if $\dim M = 4$ and 
$f:M \longrightarrow \BC^3$ is an immersion, then we have 
the formula
$$
p_1(\Omega_M) = -{\Bbb C}r (f) + dT
$$
where ${\Bbb C}r (f)$ is the (generically finite) set of 
complex
tangencies to $f(M)\subset\BC^3$ taken with appropriate 
indices.
More complicated formulas involving higher-order complex 
tangencies 
and Shur functions
are derived in \cite{HL2}.

Another application associates geometric currents to pairs 
of complex
structures. Suppose $J_1$ and $J_2$ are smooth almost 
complex
structures on a vector bundle $E\lra X$, and let
$E\ox\BC=E_1\op\ov E_1=E_2\op\ov E_2$ be the associated 
splittings.
Restriction and projection give a complex bundle map $p: 
E_1\lra E_2$
to which there are associated degeneracy currents and 
geometric
formulas. For example, let $n=\rank(E)$ and suppose
$\la\equdef\La^n_{\,\BC} p$ is atomic. Then we have the 
characteristic
current
$$
\bbc r(J_1,J_2) = \Div(\la)
$$
which is supported in the set $\{x\in X: \ker(J_1+
J_2)\neq\{0\}\}$. Let
$D_k$ be a connection on $E$ such that $D_k(J_k)=0$, and set
$e_k=c_1(D_k)$ for $k=1$, $2$. Then for any 
$\phi(t)\in\BC[t]$, there
is an $\loc$ form $\s_\phi$ on $X$ such that  
$$
\phi(e_2) - \phi(e_1) =
\f{\phi(e_2) - \phi(e_1)}{e_2-e_1} \bbc r(J_1,J_2) + 
d\s_\phi.
$$
Similar formulas, which involve the higher degeneracy
currents, are derived in \cite{HL2}. 

The theory produces formulas of the above type in the 
theory of
foliations. It is also possible, using this theory, to 
rederive and
generalize formulas of Sid Webster [W$^*$]. It should be 
remarked that
discussions with Jon Wolfson (cf\. [Wo] and [MW]) about 
Webster's
formulas served as an inspiration for this work.  

Each geometric formula discussed in \cite{HL2} carries 
with it a 
smooth 1-parameter family of analogous 
formulas coming from 
${\overset{\rightarrow}\to D}_s$.  As $s \rightarrow 0$, 
the characteristic forms
converge to the degeneracy current.

In [HL3] the authors study residue formulas associated to 
bundle maps
of general rank. These formulas simultaneously involve 
several
of the degeneracy strata of the map, each paired with its 
own
residue form. They promote the work of R\.~MacPherson 
\cite{M1,
M2} to the level of characteristic forms and currents.

The authors are indebted to Bill Fulton and John Zweck for 
important suggestions and comments related to this work.


\Refs
\widestnumber\key{BGS1}

\ref
\key AH  \by M. F. Atiyah and F. Hirzebruch
\paper The Riemann-Roch theorem for analytic embeddings
\jour Topology     \vol1      \yr1962    \pages151--166
\endref


\ref
\key BGS1  \by J.-M. Bismut, H. Gillet, and C. Soul\'e
\paper Analytic torsion and holomorphic determinant 
bundles. {\rm I, 
II, III}
\jour Comm. Math. Physics     \vol115     \yr1988    
\pages49--78, 79--126, 301--351
\endref

\ref
\key BGS2  \bysame
\paper  Bott-Chern currents and complex immersions
\jour Duke Math. J.  \vol60   \yr1990    \pages255--284
\endref

\ref
\key BGS3  \bysame
\book  Complex immersions and Arakelov geometry
\publ Grothendieck Fetschrift I,  Birkh\"auser     
\publaddr Boston
\yr1990     \pages249--331
\endref

\ref
\key BV  \by N. Berline and M. Vergne
\paper  A computation of the equivariant index of the Dirac
operator
\jour Bull. Soc. Math. France   \vol113    \yr1985  
\pages305--345
\endref

\ref
\key C1  \by S. Chern
\paper A simple intrinsic proof of the Gauss-Bonnet 
formula for
closed Riemannian manifolds
\jour Ann. of Math. (2)   \vol45     \yr1944    
\pages747--752
\endref

\ref
\key C2  \bysame
\paper On the curvature integra in a Riemannian manifold
\jour Ann. of  Math. (2)    \vol46     \yr1945    
\pages674--684
\endref

\ref
\key F  \by H. Federer
\book Geometric measure theory
\publ Springer-Verlag \publaddr New York    \yr1969     
\endref

\ref
\key Fu  \by W. Fulton
\paper Intersection theory
\inbook Ergeb. Math. Grenzgeb. (3), 
bd. 2
\publ Springer-Verlag  \publaddr  Berlin and Heidelberg   
\yr1984
\endref

\ref
\key HL1  \by F. R. Harvey and H. B. Lawson, Jr.
\paper {\rm A theory of characteristic currents associated 
with a 
singular connection}
\jour Ast\'erisque    \yr1993   \vol213\pages 1--268
\endref

\ref
\key HL2  \bysame
\book Geometric residue theorems
\bookinfo MSRI Preprint, 1993 \pages 1--53
\endref

\ref
\key HL3  \bysame
\paper A theory of characteristic currents associated with 
a 
singular connection --- Part {\rm II}
\paperinfo in preparation
\endref

\ref
\key HP  \by F. R. Harvey and J. Polking
\paper Fundamental solutions in complex analysis, Part \RM I
\jour Duke Math. J.    \vol46    \yr1979     \pages253--300
\endref

\ref
\key HS  \by F. R. Harvey and S. Semmes
\paper Zero divisors of atomic functions
\jour Ann. of Math. (2)     \yr1992   \vol135   
\pages567--600
\endref

\ref
\key M1 \by R. MacPherson
\paper Singularities of vector bundle maps
\inbook Proceedings of Liverpool Singularities Symposium, I
\bookinfo Lecture Notes in Math.,  vol. 192
\publ Springer-Verlag \publaddr New York    \yr1971    
\pages 316--318 
\endref

\ref
\key M2 \bysame
\paper Generic vector bundle maps 
\inbook  Dynamical Systems
(Proceedings of Symposium, University of Bahia, Salvador, 
1971) 
\publ  Academic Press, New York  \yr1973     \pages 
165--175 
\endref

\ref
\key MQ  \by V. Mathai and D. Quillen
\paper Superconnections, Thom classes, and equivariant 
differential
forms
\jour Topology   \vol25   \yr1986   \pages85--110
\endref

\ref
\key MW  \by M. Micallef and J. Wolfson
\paper The second variation of area of minimal surfaces in 
four-manifolds
\jour Math. Ann.  \vol295  \yr1993  \pages245--267
\endref

\ref
\key Q  \by D. Quillen
\paper Superconnections and the Chern character 
\jour Topology    \vol24     \yr1985    \pages89--95
\endref

\ref
\key R \by F. Ronga
\paper Le calcul de la classe de cohomologie
duale a $\bar{\Sigma}^k$\<
\inbook {\rm Proceedings of Liverpool 
Singularities Symposium, I}
\bookinfo  Lecture Notes in Math.,  vol. 192
\publ Springer-Verlag\publaddr New York    \yr1971    
\pages 313--315 
\endref

\ref
\key W1  \by S. Webster
\paper  Minimal surfaces in a K\"ahler surface
\jour J. Differential Geom.    \vol20    \yr1984    
\pages463--470
\endref

\ref
\key W2  \bysame
\paper The Euler and Pontrjagin numbers of an $n$-manifold 
in 
$\BC^n$
\jour Comment. Math. Helv.   \vol60    \yr1985    
\pages193--216
\endref

\ref
\key W3  \bysame
\paper On the relation between Chern and Pontrjagin numbers 
\jour Contemp. Math. \vol 49     \yr 1986    \pages135--143
\endref

\ref
\key Wo  \by J. Wolfson
\paper  Minimal surfaces in K\"ahler surfaces and Ricci 
curvature
\jour J. Differential Geom.   \vol29    \yr1989    
\pages281--294
\endref

\ref
\key Z1  \by J. Zweck
\paper Chern currents of singular connections associated 
with a
section of a compactified bundle
\paperinfo 1--40 (to appear)
\endref

\ref
\key Z2  \bysame
\paper Euler and Pontrjagin currents of a section of a 
compactified
real bundle
\paperinfo 1--32 (to appear)
\endref

\endRefs
\enddocument